\documentclass[12pt]{amsart}
\usepackage{amssymb}
\usepackage[all]{xy}
\usepackage{tikz}
\usepackage{amsmath}
\numberwithin{equation}{section}
\input xy
\xyoption{all}

\newtheorem{lemma}{Lemma}[section]
\newtheorem{corollary}[lemma]{Corollary}
\newtheorem{theorem}[lemma]{Theorem}
\newtheorem{proposition}[lemma]{Proposition}

\theoremstyle{definition}

\newtheorem{definition}[lemma]{Definition}

\DeclareMathOperator{\Mod}{Mod}
\DeclareMathOperator{\modd}{mod}

\DeclareMathOperator{\Hom}{Hom}

\DeclareMathOperator{\Ext}{Ext}

\DeclareMathOperator{\Lim}{Lim}

\DeclareMathOperator{\Ker}{Ker}
\DeclareMathOperator{\Coker}{Coker}
\DeclareMathOperator{\Imm}{Im}
\DeclareMathOperator{\fp}{fp}
\DeclareMathOperator{\Flatt}{Flat}

\DeclareMathOperator{\Eff}{Eff}
\DeclareMathOperator{\eff}{eff}

\newtheorem*{theorem a*}{Theorem A}
\newtheorem*{theorem b*}{Theorem B}

\begin{document}

\title{Yoneda extensions of abelian quotient categories}

\author{Ramin Ebrahimi} 
\address{Department of Pure Mathematics\\
Faculty of Mathematics and Statistics\\
University of Isfahan\\
P.O. Box: 81746-73441, Isfahan, Iran\\ and School of Mathematics, Institute for Research in Fundamental Sciences (IPM), P.O. Box: 19395-5746, Tehran, Iran}
\email{ramin69@sci.ui.ac.ir / ramin.ebrahimi1369@gmail.com}



\begin{abstract}
Let $\mathcal{A}$ be a essentially small abelian category and $\mathcal{C}$ be a Serre subcategory of  $\mathcal{A}$. Consider the quotient functor $q:\mathcal{A}\rightarrow {\mathcal{A}}/{\mathcal{C}}$. For an object $A\in \mathcal{A}$ and a non-negative integer $k$ we investigate when the natural map $q^i_{X,A}:\Ext_{\mathcal{A}}^i(X,A)\rightarrow\Ext_{{\mathcal{A}}/{\mathcal{C}}}^i(q(X),q(A))$ is invertible for every $X\in \mathcal{A}$ and every $i\in \{0,1,\cdots,k\}$. In the end we give an application of the main theorem.
\end{abstract}

\maketitle

\section{Introduction}
Let $\mathcal{A}$ be an abelian category. A Serre subcategory $\mathcal{C}\subseteq\mathcal{A}$ is a full subcategory that is closed under subobjects, quotient objects, and
extensions. In this case there exist an abelian category ${\mathcal{A}}/{\mathcal{C}}$ and an exacts functor
\begin{equation}
q:\mathcal{A}\rightarrow {\mathcal{A}}/{\mathcal{C}} \notag
\end{equation}
that is universal among exact functors from $\mathcal{A}$ to an abelian category that annihilate $\mathcal{C}$.
The quotient category ${\mathcal{A}}/{\mathcal{C}}$ was defined by Grothendieck and studied extensively by Gabriel \cite{G}.

The quotient functor $q$ is not full or faithful, but by a result of Gabriel for each object $A$ in the left perpendicular subcategory of $\mathcal{C}$, the natural map $q_{X,A}:\Hom_{\mathcal{A}}(X,A)\rightarrow\Hom_{{\mathcal{A}}/{\mathcal{C}}}(q(X),q(A))$ is invertible, for every $X\in \mathcal{A}$.

In many situations $\mathcal{A}$ is a known and well behaved category like a module category or a functor category and we can compute $\Ext$-groups for any two objects in $\mathcal{A}$. But the quotient category ${\mathcal{A}}/{\mathcal{C}}$ has a complicated structure, for instance if $\mathcal{C}$ is a localizing subcategory (i.e. $q$ has a right adjoint) then ${\mathcal{A}}/{\mathcal{C}}$ is equivalent to a subcategory of $\mathcal{A}$, but not an exact subcategory, because the right adjoint is fully faithful and left exact, but not exact in general \cite[Lemma 2.2.10]{KrB}. Thus it is hard to relate the $\Ext$-groups of these two categories.

The main theorem of this paper is the following result.

\begin{theorem}\label{1.1}
Let $\mathcal{A}$ be a essentially small abelian category, $\mathcal{C}$ a Serre subcategory of $\mathcal{A}$ and $q:\mathcal{A}\rightarrow {\mathcal{A}}/{\mathcal{C}}$ be the quotient functor. For an object $A\in \mathcal{A}$ and a non-negative integer $k$ the following statements are equivalent.
\begin{itemize}
\item[$(i)$]
$A\in \mathcal{C}^{{\bot}_{k+1}}$. i.e. $\Ext_{\mathcal{A}}^i(\mathcal{C},A)=0$ for every $i\in \{0,1,\cdots,k+1\}$.
\item[$(ii)$]
The natural map $q^i_{X,A}:\Ext_{\mathcal{A}}^i(X,A)\rightarrow\Ext_{{\mathcal{A}}/{\mathcal{C}}}^i(q(X),q(A))$ is invertible, for every $X\in \mathcal{A}$ and every $i\in \{0,1,\cdots,k\}$.
\end{itemize}
\end{theorem}

Note that C. Psaroudakis proves this result for a recollement of abelian categories with enough projective and injective objects \cite{Ps}. In other word we prove that the existence of left and right adjoint for the quotient functor
$q:\mathcal{A}\rightarrow {\mathcal{A}}/{\mathcal{C}}$ and also the existence of enough injective and enough projective objects are not necessary.


\section{preliminaries}
\subsection{Localization of abelian categories}
We recall below the basic theory of localisation of abelian categories. The proofs can be found in Gabriel's
thesis \cite{G} or in standard textbooks like Faith's book \cite{F} or
the recent book of Krause \cite{KrB}.

Let $\mathcal{A}$ be an abelian category. A full subcategory $\mathcal{C}$ of $\mathcal{A}$ is called a \textbf{Serre subcategory} if for
any short exact sequence
\begin{equation}
0\rightarrow A_1\rightarrow A_2\rightarrow A_3\rightarrow 0 \notag
\end{equation}
we have that $A_2\in \mathcal{C}$ if and only if $A_1, A_3\in \mathcal{C}$. In this case we have the quotient category ${\mathcal{A}}/{\mathcal{C}}$ that is by definition localization of $\mathcal{A}$ with respect to the class of all morphisms $f:X\rightarrow Y$ with $\Ker(f), \Coker(f)\in \mathcal{C}$.
The quotient category ${\mathcal{A}}/{\mathcal{C}}$ is an abelian category and there is a universal exact functor
\begin{center}
$q:\mathcal{A}\longrightarrow \mathcal{A}/\mathcal{C}$
\end{center}
such that $q(X) = 0$ if and only if $X\in \mathcal{C}$. Furthermore any other exact functor $F:\mathcal{A}\rightarrow \mathcal{D}$ annihilating $\mathcal{C}$ where $\mathcal{D}$ is an abelian category factor uniquely through $q$.

A Serre subcategory $\mathcal{C}\subseteq \mathcal{A}$ is called a \textbf{localizing subcategory} if the quotient functor $q:\mathcal{A}\rightarrow {\mathcal{A}}/{\mathcal{C}}$ admits a right adjoint
$r:{\mathcal{A}}/{\mathcal{C}}\rightarrow \mathcal{A}$. The right adjoint $r$, called the \textbf{section functor}, is fully faithful \cite[Lemma 2.2.10]{KrB}. A Serre subcategory of a Grothendieck category is localizing if and only if it is closed under coproducts \cite[Corollary 2.2.17]{KrB}. Another useful property of Grothendieck categories is the fact that every object $A$ admits an injective envelope $A\hookrightarrow I$ \cite[Corollary 2.5.4]{KrB}. This means that $A$ is a subobject of an injective object $I$, and $I$ is an essential extension of $A$ (i.e. the intersection of $A$ with every non-zero subobject of $I$ is non-zero).

Let $\mathcal{C}$ be a Serre subcategory of an abelian category $\mathcal{A}$. Recall that an object $A\in \mathcal{A}$ is called $\mathcal{C}$-\textbf{closed}
if for every morphism $f : X \rightarrow Y$ with $\Ker(f), \Coker(f)\in \mathcal{C}$ we have that
$\Hom_{\mathcal{A}}(f, A)$ is bijective. Denote by $\mathcal{C}^{\bot}$ the full subcategory of all $\mathcal{C}$-closed objects, the following result is well known. See for instance the fundamental paper of Geigle and Lenzing \cite{GL} or the recent book of Krause \cite{KrB}.

\begin{theorem}\label{2.1}
Let $\mathcal{C}$ be a Serre subcategory of an abelian category $\mathcal{A}$. The following
statements hold:
\begin{itemize}
\item[$(i)$]
We have
\begin{center}
$\mathcal{C}^{\bot}=\{A\in \mathcal{A} \mid \Hom_{\mathcal{A}}(\mathcal{C},A)=0=\Ext^1_{\mathcal{A}}(\mathcal{C},A)\}.$
\end{center}
\item[$(ii)$]
The natural map $q_{X,A}:\Hom_{\mathcal{A}}(X,A)\rightarrow \Hom_{{\mathcal{A}}/{\mathcal{C}}}(q(X),q(A))$ is invertible for every $X\in \mathcal{A}$ if and only $A\in \mathcal{C}^{\bot}$.
\item[$(iii)$]
If $\mathcal{C}$ is a localizing subcategory, the restriction
$q:\mathcal{C}^{\bot}\rightarrow {\mathcal{A}}/{\mathcal{C}}$ is an equivalence of categories. And it's quasi-inverse is induced by the section functor $r:{\mathcal{A}}/{\mathcal{C}}\rightarrow \mathcal{C}^{\bot}\subseteq \mathcal{A}$.
\item[$(iv)$]
If $\mathcal{C}$ is localizing and $\mathcal{A}$ has injective envelopes, then $\mathcal{C}^{\bot}$ has injective envelopes and the inclusion functor $\mathcal{C}^{\bot}\hookrightarrow \mathcal{A}$ preserves injective envelopes. 
\end{itemize}
\end{theorem}

The following definition is borrowed from \cite{Ps}.

\begin{definition}\label{2.2}
Let $\mathcal{A}$ be an abelian category and $\mathcal{X}$ be a subcategory of $\mathcal{A}$. For a positive integer $k$ we denote by $\mathcal{X}^{{\bot}_k}$ the full subcategory of $\mathcal{A}$ defined by $\mathcal{X}^{{\bot}_k} = \{A \in \mathcal{A} | \Ext^{0,...,k}(\mathcal{X}, A) = 0\}$. $^{{\bot}_k}\mathcal{X}$ is defined similarly \cite{Ps}. 
\end{definition}

Note that by Theorem \ref{2.1} for a Serre subcategory $\mathcal{C}$ of $\mathcal{A}$ we have $\mathcal{C}^{{\bot}_1}=\mathcal{C}^{\bot}$.

Let $\mathcal{B}$ be an abelian category and $\mathcal{S}$ be a localizing subcategory of $\mathcal{B}$. Denote the localization functor and the section functor with $e$ and $r$ respectively, we have the following diagram of functors.
\begin{center}
\begin{tikzpicture}
\node (X2) at (0,0) {$\mathcal{B}$};
\node (X3) at (3,0) {${\mathcal{B}}/{\mathcal{S}}$};
\node (r) at (1.5,-0.7) {$r$};
\draw [->,thick] (X2) -- (X3) node [midway,above] {$e$};
\draw [->,thick] (X3) to [out=225,in=315] (X2) node [midway,left] {};
\end{tikzpicture}
\end{center}

The following proposition can be seen as a variation of Proposition 3.4 and Theorem 3.10 from \cite{Ps}.
The $\Hom$ vanishing is the crucial step of the proof of Proposition 3.3 in Psaroudakis paper (the dual result - Proposition 3.4). We will prove this condition using the existence of injective envelopes in Grothendieck categories.

\begin{proposition}\label{2.3}
Let $\mathcal{B}$ be a Grothendieck category and $\mathcal{S}$ be a localizing subcategory of $\mathcal{B}$. For an object $A\in \mathcal{B}$ and a non-negative integer $k$ the following statements are equivalent.
\begin{itemize}
\item[$(i)$]
$A\in \mathcal{S}^{{\bot}_{k+1}}$.
\item[$(ii)$]
There exists an injective coresolution
\begin{equation}\label{ri}
0\rightarrow A\rightarrow r(I^0) \rightarrow r(I^1)\rightarrow \cdots\rightarrow r(I^{k+1}) 
\end{equation}
for $A$, where $I^i$'s are injective objects in ${\mathcal{B}}/{\mathcal{S}}$.
\item[$(iii)$]
The natural map $e_{X,A}^i:\Ext_{\mathcal{B}}^i(X,A)\rightarrow\Ext_{{\mathcal{B}}/{\mathcal{S}}}^i(e(X),e(A))$ is invertible, for every $X\in \mathcal{B}$ and every $i\in\{0,1,\cdots,k\}$.
\end{itemize}
\begin{proof}
$(i)\Rightarrow(ii)$
Let $A\in \mathcal{S}^{{\bot}_{k+1}}$, and $A\hookrightarrow I^0$ be the injective envelope of $A$ in $\mathcal{B}$. Let $S\in \mathcal{S}$ and $f\in \Hom_{\mathcal{B}}(S,I^0)$ be a non-zero morphism, because $I^0$ is an essential extension of $A$ we have that $\Imm(f)\cap A\neq 0$. By the definition of Serre subcategory $\Imm(f)\cap A \in \mathcal{S}$, and this contradicts the assumption
$\Hom_{\mathcal{B}}(\mathcal{S},A)=0$. Thus $I^0\in\mathcal{S}^{\bot}$ because $I^0$ is injective.

Now by applying the functor $\Hom_{\mathcal{B}}(S,-)$ for an arbitrary object $S\in \mathcal{S}$, to the short exact sequence
\begin{center}
$0\rightarrow A\rightarrow I^0\rightarrow \Omega^{-1}A\rightarrow 0$
\end{center}
we obtain the long exact sequence
\begin{align*}
0&\rightarrow \Hom_{\mathcal{B}}(S,A)\rightarrow\Hom_{\mathcal{B}}(S,I^0)\rightarrow\Hom_{\mathcal{B}}(S,\Omega^{-1}A)\\
&\rightarrow\Ext_{\mathcal{B}}^1(S,A)\rightarrow\Ext_{\mathcal{B}}^1(S,I^0)\rightarrow\Ext_{\mathcal{B}}^1(S,\Omega^{-1}A)\\
&\rightarrow\Ext_{\mathcal{B}}^2(S,A)\cdots .
\end{align*}
By assumption $\Hom_{\mathcal{B}}(S,A)=\Ext_{\mathcal{B}}^1(S,A)=\cdots=\Ext_{\mathcal{B}}^{k+1}(S,A)=0$, and $\Ext_{\mathcal{B}}^1(S,I^0)=0$, because $I^0$ is injective.
Using these vanishing conditions and the above long exact sequence we see that $\Omega^{-1}A\in\mathcal{S}^{\bot_k}$.

Now consider the injective envelope $\Omega^{-1}A\hookrightarrow I^1$ for $\Omega^{-1}A$ and the induced short exact sequence
\begin{center}
$0\rightarrow \Omega^{-1}A\rightarrow I^1\rightarrow \Omega^{-2}A\rightarrow 0.$
\end{center}
Applying the above argument for this short exact sequence we see that $I^1\in \mathcal{S}^{\bot}$ and $\Omega^{-2}A\in\mathcal{S}^{\bot_{k-1}}$. By repeating this argument and using the dimension shifting argument we obtain an injective coresolution as \eqref{ri} inductively.

$(ii)\Rightarrow(iii)$
Apply the functor $\Hom_{\mathcal{B}}(X,-)$ to the injective coresolution \eqref{ri} and use the fact that $(e,r)$ is an adjoint pair.

$(iii)\Rightarrow(i)$
By assumption for every $S\in \mathcal{S}$ and every $i\in\{0,1,\cdots,k\}$, we have
\begin{center}
$\Ext_{\mathcal{B}}^i(S,A)\cong\Ext_{{\mathcal{B}}/{\mathcal{S}}}^i(e(S),e(A))=
\Ext_{{\mathcal{B}}/{\mathcal{S}}}^i(0,e(A))=0.$
\end{center}
This means that $A\in \mathcal{S}^{\bot_k}$, so it remains to show that $\Ext_{\mathcal{B}}^{k+1}(\mathcal{S},A)=0$. We prove the later by induction. For $k=0$ the result follows from Theorem \ref{2.1}. Assume that $k\geq 1$ be a positive integer, and the claim is true for $k-1$.
Consider exact sequence 
\begin{center}
$0\rightarrow A\rightarrow I\rightarrow B\rightarrow 0$
\end{center}
Where $I$ is the injective envelope of $A$. We need to show that $\Ext_{\mathcal{B}}^{k}(\mathcal{S},B)=0$.

As we saw in above, $I\in \mathcal{S}^{\bot}$ and
$B\in \mathcal{S}^{\bot_{k-1}}$ because $A\in \mathcal{S}^{\bot_{k}}$.
Thus $\Ext_{\mathcal{B}}^{k}(\mathcal{S},B)=0$ if and only if $B\in \mathcal{S}^{\bot_{k}}$ and this is by induction hypothesis, equivalent to the condition that natural map
\begin{center}
$e_{X,B}^i:\Ext_{\mathcal{B}}^i(X,B)\rightarrow\Ext_{{\mathcal{B}}/{\mathcal{S}}}^i(e(X),e(B))$
\end{center}
be invertible, for every $X\in \mathcal{B}$ and every $i\in\{0,1,\cdots,k-1\}$.
Because $B\in \mathcal{S}^{\bot_{k-1}}$, $e_{X,B}^i$ is invertible for every $i\in\{0,1,\cdots,k-2\}$. Thus we only need to show that $e_{X,B}^{k-1}$ is invertible. Applying the functors $\Hom_{\mathcal{B}}(X,-)$ and $\Hom_{\mathcal{B}/\mathcal{S}}(e(X),-)$ to short exact sequences $0\rightarrow A\rightarrow I\rightarrow B\rightarrow 0$ and $0\rightarrow e(A)\rightarrow e(I)\rightarrow (B)\rightarrow 0$ respectively we obtain the following commutative diagram with exact rows.
\begin{center}
\begin{tikzpicture}
\node (X1) at (-8,0) {$\Ext_{\mathcal{B}}^{k-1}(X,A)$};
\node (X2) at (-4,0) {$\Ext_{\mathcal{B}}^{k-1}(X,I)$};
\node (X3) at (0,0) {$\Ext_{\mathcal{B}}^{k-1}(X,B)$};
\node (X4) at (4,0) {$\Ext_{\mathcal{B}}^{k}(X,A)$};
\node (X5) at (6.2,0) {$0$};
\node (X6) at (-8,-2) {$\Ext^{k-1}(e(X),e(A))$};
\node (X7) at (-4,-2) {$\Ext^{k-1}(e(X),e(I))$};
\node (X8) at (0,-2) {$\Ext^{k-1}(e(X),e(B))$};
\node (X9) at (4,-2) {$\Ext^{k}(e(X),e(A))$};
\node (X10) at (6.2,-2) {$0$};
\draw [->,thick] (X1) -- (X2) node [midway,left] {};
\draw [->,thick] (X2) -- (X3) node [midway,right] {};
\draw [->,thick] (X3) -- (X4) node [midway,above] {};
\draw [->,thick] (X4) -- (X5) node [midway,above] {};
\draw [->,thick] (X6) -- (X7) node [midway,left] {};
\draw [->,thick] (X7) -- (X8) node [midway,right] {};
\draw [->,thick] (X8) -- (X9) node [midway,left] {};
\draw [->,thick] (X9) -- (X10) node [midway,right] {};
\draw [->,thick] (X1) -- (X6) node [midway,left] {$e_{X,A}^{k-1}$};
\draw [->,thick] (X2) -- (X7) node [midway,left] {$e_{X,I}^{k-1}$};
\draw [->,thick] (X3) -- (X8) node [midway,left] {$e_{X,B}^{k-1}$};
\draw [->,thick] (X4) -- (X9) node [midway,left] {$e_{X,A}^{k}$};
\end{tikzpicture}
\end{center}
By assumption $e_{X,A}^{k-1}$ and $e_{X,A}^{k}$ are invertible, and because $I\in \mathcal{S}^{\bot_k}$, $e_{X,I}^{k-1}$ is also invertible. Thus the five Lemma yields that $e_{X,B}^{k-1}$ is an isomorphism.
\end{proof}
\end{proposition}

Let $\mathcal{A}$ be a essentially small additive category. Recall that $\Mod\mathcal{A}$ is
the category of all additive contravariant functors from $\mathcal{A}$ to the category of all abelian
groups. It is an abelian category with all limits and colimits, which are defined point-wise. A functor $F\in \Mod\mathcal{A}$ is called finitely presented (or coherent) if there exists an exact sequence
\begin{center}
$\Hom(-,X)\rightarrow \Hom(-,Y)\rightarrow F\rightarrow 0$
\end{center}
in $\Mod\mathcal{A}$ \cite{Au}. We denote by $\modd\mathcal{A}$ the full subcategory of $\Mod\mathcal{A}$ consists of all finitely presented functors. It is a well known result that $\modd\mathcal{A}$ is an abelian category and the inclusion $\modd\mathcal{A}\hookrightarrow \Mod\mathcal{A}$ is an exact functor if and only if $\mathcal{A}$ have weak kernels \cite[Lemma 2.1.6]{KrB}. In particular this is the case when $\mathcal{A}$ is an abelian category.
By the Yoneda lemma, representable functors are projective and the direct sum of all representable functors $\bigoplus_{X\in \mathcal{M}}\Hom(-,X)$, is a generator for $\Mod\mathcal{A}$. Thus $\Mod\mathcal{A}$ is a Grothendieck category \cite[Proposition 5.21]{Fr}.

\begin{definition}\label{2.4}
\begin{itemize}
\item[$(i)$]
A functor $F\in \Mod\mathcal{A}$ is called \textbf{weakly effaceable}, if for each object $X\in \mathcal{A}$ and $x\in F(X)$
there exists an epimorphism $f : Y \rightarrow X$ such that $F(f)(x) = 0$ \cite[Page 28]{We}. We denote by
$\Eff(\mathcal{A})$ the full subcategory of all weakly effaceable functors.
\item[$(ii)$]
A functor $F\in \modd\mathcal{A}$ is called \textbf{effaceable}, if there exists an exact sequence
\begin{center}
$\Hom(-,Y)\rightarrow \Hom(-,X)\rightarrow F\rightarrow 0$
\end{center}
such that $Y \rightarrow X$ is an epimorphism. We denote by
$\eff(\mathcal{A})$ the full subcategory of all effaceable functors.
\item[$(iii)$]
A functor $F\in \Mod\mathcal{A}$ is called a left exact functor, if for each short exact sequence $0\rightarrow X\rightarrow Y\rightarrow Z\rightarrow 0$ in $\mathcal{A}$, the sequence of abelian groups $0\rightarrow F(Z)\rightarrow F(Y)\rightarrow F(X)$ is exact. We denote by $\mathcal{L}(\mathcal{A})$ the full subcategory of all left exact functors.
\end{itemize}
\end{definition}
Note that Krause in his recent book \cite{KrB}, has used the terminology \textbf{locally effaceable} instead of weakly effaceable. And this is reasonable from the viewpoint of locally finitely presented categories (see the next subsection), because we can show that objects in $\Eff(\mathcal{A})$ are exactly direct limits of objects in $\eff(\mathcal{A})$ \cite[Page 672]{Kr15}.

\begin{proposition}\label{2.5}
Let $\mathcal{A}$ be a essentially small abelian category.
\begin{itemize}
\item[$(i)$]
$\rm Eff(\mathcal{A})$ is a localizing subcategory of $\rm{Mod}\mathcal{A}$.
\item[$(ii)$]
We have 
\begin{align*}
\mathcal{L}(\mathcal{A})&=\Eff(\mathcal{A})^{\bot}\simeq (\Mod\mathcal{A})/(\Eff(\mathcal{A})),
\end{align*}
and the inclusion $\mathcal{L}(\mathcal{A})\hookrightarrow \Mod\mathcal{A}$ is the right adjoint of the localization functor $e:\Mod\mathcal{A}\rightarrow \mathcal{L}(\mathcal{A})$.
\item[$(iii)$]
$e(G)=0$ if and only if $G\in \Eff(\mathcal{A})$.
\end{itemize}
\begin{proof}
This is a standard result proved by Gabriel \cite[II.2]{G} (see also \cite[Proposition 2.3.7]{KrB}).
\end{proof}
\end{proposition}

We denote by $\mathsf{i}:\mathcal{A}\rightarrow\mathcal{L}(\mathcal{A})$ the composition of the Yoneda functor $\mathcal{A}\rightarrow \Mod \mathcal{A}$ and the localization functor $\Mod \mathcal{A}\rightarrow (\Mod \mathcal{A}/(\Eff(\mathcal{A}))\simeq \Eff(\mathcal{A})^{\bot}=\mathcal{L}(\mathcal{A})$.
The canonical functor $\mathsf{i}$ is exact \cite[Proposition 2.3.7]{KrB}.
For an object $A\in \mathcal{A}$ sometimes we denote $\mathsf{i}(A)=\Hom_{\mathcal{A}}(-,A)$ by $H_A$.
The following easy lemma will be used throughout the paper.

\begin{lemma}\label{2.6}
Let $A\in \mathcal{A}$ and $\alpha:F\rightarrow H_A$ be an epimorphism in $\mathcal{L}(\mathcal{A})$. Then the cokernel of $\alpha$ computed in $\Mod\mathcal{A}$ is weakly effaceable. In particular, there exists an epimorphism $B\rightarrow A$ inducing the following commutative diagram with exact rows in $\mathcal{L}(\mathcal{A})$.
\begin{center}
\begin{tikzpicture}
\node (X0) at (-2,2) {$H_B$};
\node (X1) at (0,2) {$H_A$};
\node (X2) at (-2,0) {$F$};
\node (X3) at (0,0) {$H_A$};
\node (X4) at (2,2) {$0$};
\node (X5) at (2,0) {$0$};
\draw [->,thick] (X0) -- (X1) node [midway,left] {};
\draw [->,thick] (X0) -- (X2) node [midway,right] {};
\draw [double,-,thick] (X1) -- (X3) node [midway,above] {};
\draw [->,thick] (X2) -- (X3) node [midway,above] {};
\draw [->,thick] (X1) -- (X4) node [midway,left] {};
\draw [->,thick] (X3) -- (X5) node [midway,right] {};
\end{tikzpicture}
\end{center}
\begin{proof}
Denote by $G$ the cokernel of $\alpha$ in $\Mod\mathcal{A}$. So we have the exact sequence $F\rightarrow H_A\rightarrow G\rightarrow 0$ in $\Mod\mathcal{A}$. By applying the localization functor $e:\Mod\mathcal{A}\rightarrow \mathcal{L}(\mathcal{A})$ to this exact sequence and using Proposition \ref{2.5} we conclude that $G$ is weakly effaceable. Now consider the image of $1_A$, the identity morphism on $A$, under the map $\Hom_{\mathcal{A}}(A,A)\rightarrow G(A)$ and denote it by $a$. Since $G$ is weakly effaceable, by definition there exists an epimorphism $f:B\rightarrow A$ such that $G(f)(a)=0$. This means that $H_B\rightarrow H_A$ factor through $\alpha$, and so we have the desire diagram (cf. \cite[Page 409]{Kl}).
\end{proof}
\end{lemma}

\subsection{Locally coherent categories}
In this subsection we recall some of the basic properties of Locally finitely presented categories, for details and more information the reader is referred to \cite{CB}.
Recall that a non-empty
category $\mathcal{I}$ is said to be filtered provided that for each pair of objects $\lambda_1,\lambda_2\in \mathcal{I}$ there are morphisms $\varphi_i:\lambda_i\rightarrow \mu$ for some $\mu \in \mathcal{I}$, and for each pair of morphisms $\varphi_1,\varphi_2:\lambda\rightarrow \mu$ there is a morphism $\psi :\mu\rightarrow \nu$ with $\psi \varphi_1=\psi \varphi_2$. Let $\mathcal{A}$ be an additive category, $\mathcal{I}$ a small filtered category and $X:\mathcal{I}\rightarrow \mathcal{A}$ be an additive functor. As usual we use the term \textbf{direct limit} for colimit of $X$ when $I$ is filtered. And we denote it by $\underrightarrow{\Lim}_{i\in \mathcal{I}}X_i$.

\begin{definition}\label{2.7}
Let $\mathcal{B}$ be an additive category with direct limits .
\begin{itemize}
\item[$(i)$]
An object $X\in \mathcal{B}$ is called
finitely presented (finitely generated) provided that for every direct limit $\underrightarrow{\Lim}_{i\in \mathcal{I}}X_i$ in
$\mathcal{B}$ the natural morphism 
\begin{center}
$\underrightarrow{\Lim}\Hom_{\mathcal{B}}(X,Y_i)\rightarrow \Hom_{\mathcal{B}}(X,\underrightarrow{\Lim}Y_i)$
\end{center}
is an isomorphism (a monomorphism).
The full subcategory of finitely presented objects of $\mathcal{A}$ is denoted by $\fp(\mathcal{A})$. 
\item[$(ii)$]
$\mathcal{B}$ is called a locally finitely presented category if
$\fp(\mathcal{B})$ is essentially small and every object in $\mathcal{B}$ is a direct limit of objects in $\fp(\mathcal{B})$. 
\item[$(iii)$]
Let $\mathcal{B}$ be a locally finitely presented abelian category. $\mathcal{B}$ is said to be locally coherent provided that finitely generated subobjects of finitely presented objects are finitely presented.
\end{itemize}
\end{definition}
We recall a concrete description of locally finitely presented additive categories due to Crawley-Boevey \cite{CB}. Let $\mathcal{C}$ be a essentially small additive category, recall that there is a tensor product bifunctor
\begin{align}
\Mod\mathcal{C}\otimes_{\mathcal{C}}& \Mod\mathcal{C}^{op}\longrightarrow \rm Ab. \notag \\
(F,G)&\longmapsto F\otimes_{\mathcal{C}}G \notag
\end{align}
A functor $F\in \Mod\mathcal{C}$ is said to be flat provided that $F\otimes_{\mathcal{C}}-$ is an exact functor. We denote by $\Flatt(\mathcal{C})$ the full subcategory of $\Mod\mathcal{C}$ consist of flat functors. A theorem of Lazard state that an $R$-module is flat if and only if it is a direct limit of finitely generated free modules. Lazard theorem has been generalized to functors by Oberst and Rohrl \cite{OR}. Indeed a functor $F\in \Mod\mathcal{C}$ is flat if and only if $F$ is a direct limit of representable functors. Using this description of flat functors we have:

\begin{lemma}\label{2.8}
Let $\mathcal{A}$ be a essentially small abelian category and $F\in \Mod \mathcal{A}$. The following statements are equivalent.
\begin{itemize}
\item[(a)]
$F$ is a flat functor.
\item[(b)]
$F$ is a left exact functor.
\end{itemize}
In other words $\Flatt(\mathcal{A})=\mathcal{L}(\mathcal{A})$.
\begin{proof}
See \cite[Lemma 11.1.14]{KrB}.
\end{proof}
\end{lemma}

\begin{theorem}$($\cite[Theorem 1.4]{CB}$)$\label{2.9}
\begin{itemize}
\item[(a)]
If $\mathcal{C}$ is a essentially small additive category, then $\Flatt(\mathcal{C})$ is a locally finitely presented category, and $\fp(\Flatt(\mathcal{C}))$ consists of the direct summands of representable functors. If $\mathcal{C}$ has split idempotents then the Yoneda functor $h:\mathcal{C}\rightarrow \fp(\Flatt(\mathcal{C}))$ is an equivalence. 
\item[(b)]
If $\mathcal{B}$ is a locally finitely presented category then $\fp(\mathcal{B})$ is a
essentially small additive category with split idempotents, and the functor
\begin{align}
g:\mathcal{B}&\longrightarrow \Flatt(\fp(\mathcal{B}))\notag \\
M&\longmapsto \Hom(-,M)|_{\fp(\mathcal{B})}\notag 
\end{align}
is an equivalence.
\end{itemize}
\end{theorem}

\begin{proposition}\label{2.10}
\begin{itemize}
\item[$(i)$]
Every locally finitely presented abelian category is a Grothendieck category.
\item[$(ii)$]
A locally finitely presented abelian category $\mathcal{B}$ is locally coherent if and only if
$\fp(\mathcal{B})$ is an abelian category.
\end{itemize}
\begin{proof}
We refer the reader to \cite{CB} for the statement $(i)$ and to \cite{Ro} for $(ii)$.
\end{proof}
\end{proposition}
In the proof of the main theorem we need the following proposition.


\section{Proof of the main theorem}
Let $\mathcal{A}$ be a essentially small abelian category and $\mathcal{C}$ be a Serre subcategory of $\mathcal{A}$.
By Theorem \ref{2.9}, $\mathcal{L}(\mathcal{A})\simeq(\Mod \mathcal{A})/(\Eff(\mathcal{A}))$ is a locally coherent category and the essential image of the canoonical functor
\begin{equation}
\mathsf{i}:\mathcal{A}\longrightarrow \mathcal{L}(\mathcal{A}) \notag
\end{equation}
denoted by $\mathsf{i}(\mathcal{A})$ is the subcategory of finitely presented objects. Thus by \cite[Theorem 2.8]{Kr97} $\overrightarrow{\mathsf{i}(\mathcal{C})}$, the subcategory of $\mathcal{L}(\mathcal{A})$ consists of all direct limits of objects in $\mathsf{i}(\mathcal{C})$ is a localizing subcategory of $\mathcal{L}(\mathcal{A})$. Also by \cite[Proposition A5]{Kr01} $(\mathcal{L}(\mathcal{A}))/(\overrightarrow{{\mathsf{i}(\mathcal{C})}})$ is a locally coherent category and it's subcategory of finitely presented objects is equivalent to $\mathcal{A}/\mathcal{C}$.
Thus we have the commutative diagram of functors
\begin{equation}
\begin{tikzpicture}\label{D}
\node (X0) at (-1,3) {$\mathcal{A}$};
\node (X1) at (4,3) {$\mathcal{A}/\mathcal{C}$};
\node (X2) at (-1,0) {$\mathcal{L}(\mathcal{A})$};
\node (X3) at (4,0) {$(\mathcal{L}(\mathcal{A}))/(\overrightarrow{\mathsf{i}(\mathcal{C})})\simeq \mathcal{L}(\mathcal{A}/\mathcal{C})$};
\node (r) at (1.5,-1.2) {$r$};
\draw [->,thick] (X0) -- (X2) node [midway,left] {$\mathsf{i}$};
\draw [->,thick] (X1) -- (X3) node [midway,right] {$\mathsf{j}$};
\draw [->,thick] (X0) -- (X1) node [midway,above] {$q$};
\draw [->,thick] (X2) -- (X3) node [midway,above] {$e$};
\draw [->,thick] (X3) to [out=225,in=315] (X2) node [midway,left] {};
\end{tikzpicture}
\end{equation}
satisfies
\begin{itemize}
\item[$(i)$]
$q$ and $e$ are quotient functors.
\item[$(ii)$]
$(e,r)$ is an adjoint pair and $r$ is fully faithful.
\item[$(iii)$]
Both of vertical functors $\mathsf{i}$ and $\mathsf{j}$ are the canonical functors from an abelian category to the abelian category of left exact functors.
\end{itemize}
The following lemma due to Mitchell \cite{Mi68} is one of the crucial key steps for proving the main result of the paper. For the convenient of the reader we write a self contained proof.

For more information about Yoneda extension groups we refer the reader to \cite[Chapter VII]{Mi}. We just recall that for two objects $A, B\in \mathcal{A}$, a $k$-fold extension of $A$ by $B$ is an exact sequence of the form
\begin{equation}
\xi:0\rightarrow B\rightarrow X^k\rightarrow X^{k-1}\rightarrow \cdots\rightarrow X^1\rightarrow A\rightarrow 0, \notag
\end{equation}
and $\Ext_{\mathcal{A}}^k(A,B)$ is the set of all Yoneda equivalence classes of $k$-fold extensions of $A$ by $B$.
For a morphism $f:B\rightarrow B'$, by taking push out along $f$ we obtain a $k$-fold extension of $A$ by $B'$ that we denote by $f\xi$. In a similar way and using pull back we can define $\xi g$ for a morphism $g:A'\rightarrow A$. Now assume that we have a $k$-fold extension $\xi$ of $A$ by $B$, and an $n$-fold extension $\eta$ of $B$ by $C$. Then we can splice them together and obtain a $(k+m)$-fold extension of $A$ by $C$ that we denote by $\eta\circ \xi$.

\begin{lemma}\label{3.1} 
Let $\mathcal{A}$ be a essentially small abelian category. The canonical functor $\mathsf{i}:\mathcal{A}\longrightarrow \mathcal{L}(\mathcal{A})$ is an $\Ext$-preserving functor. i.e. for every two object $A,B\in \mathcal{A}$ and every non-negative integer $i$ the natural map $\mathsf{i}_{A,B}^i:\Ext_{\mathcal{A}}^i(A,B)\rightarrow \Ext_{\mathcal{L}(\mathcal{A})}^i(\mathsf{i}(A),\mathsf{i}(B))$ is invertible.
\begin{proof}
By the Yoneda lemma $\mathsf{i}$ is fully faithful, so for $i=0$ the claim follows. Now let $i\geq 1$ be a positive integer. First we prove by induction that $\mathsf{i}_{A,B}^i$ is surjective. Consider an element $\xi \in \Ext_{\mathcal{L}(\mathcal{A})}^i(\mathsf{i}(A),\mathsf{i}(B))$ represented as
\begin{equation}
0\rightarrow \mathsf{i}(B)\rightarrow F^i\rightarrow F^{i-1}\rightarrow \cdots\rightarrow F^1\rightarrow \mathsf{i}(A)\rightarrow 0. \notag
\end{equation}
Decompose $\xi$ as splicing of extensions
\begin{equation}
\xi_{i-1}:0\rightarrow \mathsf{i}(B)\rightarrow F^i\rightarrow F^{i-1}\rightarrow \cdots\rightarrow F^2\rightarrow G^2\rightarrow 0 \notag
\end{equation}
and
\begin{equation}
\xi_1:0\rightarrow G^2\rightarrow F^1\rightarrow \mathsf{i}(A)\rightarrow 0 . \notag
\end{equation}
By Lemma \ref{2.6} there is an exact sequence $0\rightarrow Y^2\rightarrow X^1\rightarrow A\rightarrow 0$ inducing the following commutative diagram with exact rows ($\varphi$ is induced by the universal property of kernel).
\begin{center}
\begin{tikzpicture}
\node (X0) at (-4.5,1.5) {$\xi_1^{\prime}:$};
\node (X1) at (-3,1.5) {$0$};
\node (X2) at (-1.5,1.5) {$\mathsf{i}(Y^2)$};
\node (X3) at (0,1.5) {$\mathsf{i}(X^1)$};
\node (X4) at (1.5,1.5) {$\mathsf{i}(A)$};
\node (X5) at (3,1.5) {$0$};
\node (X01) at (-4.5,0) {$\xi_1:$};
\node (X6) at (-3,0) {$0$};
\node (X7) at (-1.5,0) {$G^2$};
\node (X8) at (0,0) {$F^1$};
\node (X9) at (1.5,0) {$\mathsf{i}(A)$};
\node (X10) at (3,0) {$0.$};
\draw [->,thick] (X1) -- (X2) node [midway,left] {};
\draw [->,thick] (X2) -- (X3) node [midway,right] {};
\draw [->,thick] (X3) -- (X4) node [midway,left] {};
\draw [->,thick] (X4) -- (X5) node [midway,right] {};
\draw [->,thick] (X6) -- (X7) node [midway,left] {};
\draw [->,thick] (X7) -- (X8) node [midway,right] {};
\draw [->,thick] (X8) -- (X9) node [midway,left] {};
\draw [->,thick] (X9) -- (X10) node [midway,right] {};
\draw [->,thick] (X2) -- (X7) node [midway,left] {$\varphi$};
\draw [->,thick] (X3) -- (X8) node [midway,right] {};
\draw [double,-,thick] (X4) -- (X9) node [midway,above] {};
\end{tikzpicture}
\end{center}
So $\xi_1=\varphi\xi_1^{\prime}$. If $i=1$ the above diagram is of the form
\begin{center}
\begin{tikzpicture}
\node (X1) at (-3,1.5) {$0$};
\node (X2) at (-1.5,1.5) {$\mathsf{i}(Y^2)$};
\node (X3) at (0,1.5) {$\mathsf{i}(X^1)$};
\node (X4) at (1.5,1.5) {$\mathsf{i}(A)$};
\node (X5) at (3,1.5) {$0$};
\node (X6) at (-3,0) {$0$};
\node (X7) at (-1.5,0) {$\mathsf{i}(B)$};
\node (X8) at (0,0) {$F^1$};
\node (X9) at (1.5,0) {$\mathsf{i}(A)$};
\node (X10) at (3,0) {$0.$};
\draw [->,thick] (X1) -- (X2) node [midway,left] {};
\draw [->,thick] (X2) -- (X3) node [midway,right] {};
\draw [->,thick] (X3) -- (X4) node [midway,left] {};
\draw [->,thick] (X4) -- (X5) node [midway,right] {};
\draw [->,thick] (X6) -- (X7) node [midway,left] {};
\draw [->,thick] (X7) -- (X8) node [midway,right] {};
\draw [->,thick] (X8) -- (X9) node [midway,left] {};
\draw [->,thick] (X9) -- (X10) node [midway,right] {};
\draw [->,thick] (X2) -- (X7) node [midway,left] {$\varphi$};
\draw [->,thick] (X3) -- (X8) node [midway,right] {};
\draw [double,-,thick] (X4) -- (X9) node [midway,above] {};
\end{tikzpicture}
\end{center}
Since $\mathsf{i}$ is fully faithful, $\varphi=\mathsf{i}(f)$ for some $f:Y^2\rightarrow B$. Taking push out along $f$ and then applying the exact functor $\mathsf{i}$ we obtain the following commutative diagram with exact rows
\begin{center}
\begin{tikzpicture}
\node (X1) at (-3,1.5) {$0$};
\node (X2) at (-1.5,1.5) {$\mathsf{i}(Y^2)$};
\node (X3) at (0,1.5) {$\mathsf{i}(X^1)$};
\node (X4) at (1.5,1.5) {$\mathsf{i}(A)$};
\node (X5) at (3,1.5) {$0$};
\node (X6) at (-3,0) {$0$};
\node (X7) at (-1.5,0) {$\mathsf{i}(B)$};
\node (X8) at (0,0) {$\mathsf{i}(C)$};
\node (X9) at (1.5,0) {$\mathsf{i}(A)$};
\node (X10) at (3,0) {$0.$};
\draw [->,thick] (X1) -- (X2) node [midway,left] {};
\draw [->,thick] (X2) -- (X3) node [midway,right] {};
\draw [->,thick] (X3) -- (X4) node [midway,left] {};
\draw [->,thick] (X4) -- (X5) node [midway,right] {};
\draw [->,thick] (X6) -- (X7) node [midway,left] {};
\draw [->,thick] (X7) -- (X8) node [midway,right] {};
\draw [->,thick] (X8) -- (X9) node [midway,left] {};
\draw [->,thick] (X9) -- (X10) node [midway,right] {};
\draw [->,thick] (X2) -- (X7) node [midway,left] {$\varphi$};
\draw [->,thick] (X3) -- (X8) node [midway,right] {};
\draw [double,-,thick] (X4) -- (X9) node [midway,above] {};
\end{tikzpicture}
\end{center}
In both of above diagrams the left-hand squares are push out diagrams. Thus by the uniqueness of push out we see that short exact sequences $0\rightarrow\mathsf{i}(B)\rightarrow F^1\rightarrow\mathsf{i}(A)\rightarrow 0$ and $0\rightarrow\mathsf{i}(B)\rightarrow \mathsf{i}(C)\rightarrow\mathsf{i}(A)\rightarrow 0$ are Yoneda equivalent.

For $i\geq 2$ we have
$\xi=\xi_{i-1}\circ \xi_1=\xi_{i-1}\circ (\varphi\xi_1^{\prime})=(\xi_{i-1}\varphi)\circ\xi_1^{\prime}$.
Now by induction hypothesis $\xi_{i-1}\varphi \in \Ext_{\mathcal{L}(\mathcal{A})}^{i-1}(\mathsf{i}(Y^2),\mathsf{i}(B))$ belongs to the image of $\mathsf{i}^{i-1}_{Y^2,B}$ and $\xi_1^{\prime}\in \Ext_{\mathcal{L}(\mathcal{A})}^1(\mathsf{i}(A),\mathsf{i}(Y^2))$ belongs to the image of $\mathsf{i}^1_{A,Y^2}$. And this complete the proof of the induction step. In other word we proved that there exists a commutative diagram with exact rows of the following form.
\begin{equation}
\begin{tikzpicture}\label{DD}
\node (X3) at (-5,1) {$0$};
\node (X4) at (-3.5,1) {$\mathsf{i}(B)$};
\node (X5) at (-1.75,1) {$\mathsf{i}(X^i)$};
\node (X6) at (0,1) {$\cdots$};
\node (X7) at (1.5,1) {$\mathsf{i}(X^1)$};
\node (X8) at (3,1) {$\mathsf{i}(A)$};
\node (X9) at (4.5,1) {$0$};
\node (X10) at (-5,-0.5) {$0$};
\node (X11) at (-3.5,-0.5) {$\mathsf{i}(B)$};
\node (X12) at (-1.75,-0.5) {$F^i$};
\node (X13) at (0,-0.5) {$\cdots$};
\node (X14) at (1.5,-0.5) {$F^1$};
\node (X15) at (3,-0.5) {$\mathsf{i}(A)$};
\node (X16) at (4.5,-0.5) {$0.$};
\draw [->,thick] (X3) -- (X4) node [midway,left] {};
\draw [->,thick] (X4) -- (X5) node [midway,above] {};
\draw [->,thick] (X5) -- (X6) node [midway,above] {};
\draw [->,thick] (X6) -- (X7) node [midway,above] {};
\draw [->,thick] (X7) -- (X8) node [midway,above] {};
\draw [->,thick] (X8) -- (X9) node [midway,above] {};
\draw [->,thick] (X10) -- (X11) node [midway,above] {};
\draw [->,thick] (X11) -- (X12) node [midway,above] {};
\draw [->,thick] (X12) -- (X13) node [midway,above] {};
\draw [->,thick] (X13) -- (X14) node [midway,above] {};
\draw [->,thick] (X14) -- (X15) node [midway,above] {};
\draw [->,thick] (X15) -- (X16) node [midway,above] {};
\draw [double,-,thick] (X4) -- (X11) node [midway,above] {};
\draw [->,thick] (X5) -- (X12) node [midway,above] {};
\draw [->,thick] (X7) -- (X14) node [midway,above] {};
\draw [double,-,thick] (X8) -- (X15) node [midway,above] {};
\end{tikzpicture}
\end{equation}
Now we prove that $\mathsf{i}_{A,B}^i$ is injective. Let $\eta \in \Ext^i_{\mathcal{A}}(X,Y)$ such that $\mathsf{i}(\eta)=0$. By \cite[Theorem 4.2]{Mi} we can find a commutative diagram with exact rows
\begin{center} 
\begin{tikzpicture}
\node (xi) at (-6.5,1) {$\xi:$};
\node (eta) at (-6.5,-0.5) {$\mathsf{i}(\eta):$};
\node (X3) at (-5,1) {$0$};
\node (X4) at (-3.5,1) {$\mathsf{i}(B)$};
\node (X5) at (-1.75,1) {$F^i$};
\node (X6) at (0,1) {$\cdots$};
\node (X7) at (1.5,1) {$F^1$};
\node (X8) at (3,1) {$\mathsf{i}(A)$};
\node (X9) at (4.5,1) {$0$};
\node (X10) at (-5,-0.5) {$0$};
\node (X11) at (-3.5,-0.5) {$\mathsf{i}(B)$};
\node (X12) at (-1.75,-0.5) {$\mathsf{i}(Y^i)$};
\node (X13) at (0,-0.5) {$\cdots$};
\node (X14) at (1.5,-0.5) {$\mathsf{i}(Y^1)$};
\node (X15) at (3,-0.5) {$\mathsf{i}(A)$};
\node (X16) at (4.5,-0.5) {$0$};
\draw [->,thick] (X3) -- (X4) node [midway,left] {};
\draw [->,thick] (X4) -- (X5) node [midway,above] {};
\draw [->,thick] (X5) -- (X6) node [midway,above] {};
\draw [->,thick] (X6) -- (X7) node [midway,above] {};
\draw [->,thick] (X7) -- (X8) node [midway,above] {};
\draw [->,thick] (X8) -- (X9) node [midway,above] {};
\draw [->,thick] (X10) -- (X11) node [midway,above] {};
\draw [->,thick] (X11) -- (X12) node [midway,above] {};
\draw [->,thick] (X12) -- (X13) node [midway,above] {};
\draw [->,thick] (X13) -- (X14) node [midway,above] {};
\draw [->,thick] (X14) -- (X15) node [midway,above] {};
\draw [->,thick] (X15) -- (X16) node [midway,above] {};
\draw [double,-,thick] (X4) -- (X11) node [midway,above] {};
\draw [->,thick] (X5) -- (X12) node [midway,above] {};
\draw [->,thick] (X7) -- (X14) node [midway,above] {};
\draw [double,-,thick] (X8) -- (X15) node [midway,above] {};
\end{tikzpicture}
\end{center}
where $\mathsf{i}(B)\rightarrow F^i$ is a split monomorphism. Forming the Diagram \eqref{DD} for $\xi$ we see that $\mathsf{i}(B)\rightarrow \mathsf{i}(X^i)$ is also a split monomorphism, and so is $B\rightarrow X^i$ because $\mathsf{i}$ is fully faithful. Again because $\mathsf{i}$ is fully faithful we conclude that $\eta =0$.
\end{proof}
\end{lemma}

\begin{lemma}\label{3.2}
Let $\mathcal{B}$ be a locally coherent category and $\mathcal{C}$ be a Serre subcategory of $\mathcal{A}=\fp(\mathcal{B})$. Then $\overrightarrow{\mathcal{C}}$, the full subcategory of $\mathcal{B}$ consists of all direct limits of object in $\mathcal{A}$, is a localizing subcategory of
$\mathcal{B}$ and as subcategories of $\mathcal{B}$ we have
$\mathcal{C}^{\bot}=(\overrightarrow{\mathcal{C}})^{\bot}$.
\begin{proof}
See \cite[Theorem 2.8 and Corollary 2.11]{Kr97}.
\end{proof}
\end{lemma}

The following proposition is a generalization of the above result for arbitrary
positive integer k (Krause in \cite{Kr97} proved this result for k=1).

\begin{proposition}\label{3.3}
Let $\mathcal{B}$ be a locally coherent category, $\mathcal{C}$ be a Serre subcategory of $\mathcal{A}=\fp(\mathcal{B})$ and $k$ be a positive integer. Then as subcategories of $\mathcal{B}$ we have
$\mathcal{C}^{\bot_k}=(\overrightarrow{\mathcal{C}})^{\bot_k}$.
\begin{proof}
By Lemma \ref{3.2} we have $\mathcal{C}^{\bot_1}=(\overrightarrow{\mathcal{C}})^{\bot_1}$. Because $\mathcal{C}\subseteq \overrightarrow{\mathcal{C}}$ it is clear that $(\overrightarrow{\mathcal{C}})^{\bot_k}\subseteq \mathcal{C}^{\bot_k}$. For the converse inclusion let $X\in \mathcal{C}^{\bot_k}$. We most show that $X\in (\overrightarrow{\mathcal{C}})^{\bot_k}$. Because by Lemma \ref{3.2} $\overrightarrow{\mathcal{C}}$ is a localizing subcategory of $\mathcal{B}$, by Proposition \ref{2.3} it is enough to show that there is an injective coresolution
\begin{equation}
0\rightarrow X\rightarrow I^0\rightarrow I^1\rightarrow \cdots\rightarrow I^k\notag
\end{equation}
such that for every $ i\in\{0,1,\cdots,k\}$ we have $I^i\in (\overrightarrow{\mathcal{C}})^{\bot_1}=\mathcal{C}^{\bot_1}$.
Let $X\hookrightarrow I^0$ be the injective envelope of $X$ in $\mathcal{B}$.
Let $C\in \mathcal{C}$ and $f\in \Hom_{\mathcal{B}}(C,I^0)$ be a non-zero morphism, because $I^0$ is an essential extension of $X$ we have that $\Imm(f)\cap X\neq 0$. By the definition of Serre subcategory $\Imm(f)\cap X \in \mathcal{C}$, and this contradicts the assumption
$\Hom_{\mathcal{B}}(\mathcal{C},X)=0$. Thus $I^0\in\mathcal{C}^{\bot}$ because $I^0$ is injective.

Now by applying the functor $\Hom_{\mathcal{B}}(C,-)$ for an arbitrary object $C\in \mathcal{C}$, to the short exact sequence
\begin{center}
$0\rightarrow X\rightarrow I^0\rightarrow \Omega^{-1}X\rightarrow 0$
\end{center}
we obtain the long exact sequence
\begin{align*}
0&\rightarrow \Hom_{\mathcal{B}}(C,X)\rightarrow\Hom_{\mathcal{B}}(C,I^0)\rightarrow\Hom_{\mathcal{B}}(C,\Omega^{-1}X)\\
&\rightarrow\Ext_{\mathcal{B}}^1(C,X)\rightarrow\Ext_{\mathcal{B}}^1(C,I^0)\rightarrow\Ext_{\mathcal{B}}^1(C,\Omega^{-1}X)\\
&\rightarrow\Ext_{\mathcal{B}}^2(S,X)\cdots.
\end{align*}
By assumption $\Hom_{\mathcal{B}}(C,X)=\Ext_{\mathcal{B}}^1(C,X)=\cdots =\Ext_{\mathcal{B}}^k(C,X)=0$, and $\Ext_{\mathcal{B}}^1(C,I^0)=0$, because $I^0$ is injective.
Using these vanishing conditions and the above long exact sequence we see that $\Omega^{-1}X\in\mathcal{C}^{\bot_{k-1}}$.

Now consider the injective envelope $\Omega^{-1}X\hookrightarrow I^1$ for $\Omega^{-1}X$ and the induced short exact sequence
\begin{center}
$0\rightarrow \Omega^{-1}X\rightarrow I^1\rightarrow \Omega^{-2}X\rightarrow 0.$
\end{center}
Applying the above argument for this short exact sequence we see that $I^1\in \mathcal{C}^{\bot}$ and $\Omega^{-2}X\in\mathcal{C}^{\bot_{k-2}}$. By repeating this argument and using the dimension shifting argument we obtain the desire injective coresolution inductively.
\end{proof}
\end{proposition}

Now we are ready to prove Theorem \ref{1.1}.

\textbf{Proof of Theorem \ref{1.1}:} $(i)\Rightarrow (ii)$ Let $A\in \mathcal{C}^{{\bot}_{k+1}}$ and $X$ be an arbitrary object in $\mathcal{A}$. Having the Diagram \eqref{D} in mind by Proposition \ref{3.3}, Proposition \ref{2.3} and Lemma \ref{3.1} for every $i\in\{0,1,\cdots,k\}$ we have
\begin{align*}
\Ext_{\mathcal{A}}^i(X,A)&\overset{\mathsf{i}^i_{X,A}}{\cong} \Ext_{\mathcal{L}(\mathcal{A})}^i(\mathsf{i}(X),\mathsf{i}(A)) \\
& \overset{e^i_{X,A}}{\cong} \Ext_{\frac{\mathcal{L}(\mathcal{A})}{\overrightarrow{\mathcal{C}}}}^i(e\mathsf{i}(X),e\mathsf{i}(A)).
\end{align*}
Because $\mathsf{j}^i_{X,A} \circ q^i_{X,A}=e^i_{X,A} \circ \mathsf{i}^i_{X,A}$ we have that $\mathsf{j}^i_{X,A} \circ q^i_{X,A}$ is invertible, and since by Lemma \ref{3.1}, $\mathsf{j}^i_{X,A}$ is invertible, $q^i_{X,A}$ is also invertible.

$(ii)\Rightarrow (i)$ Assume that $q^i_{X,A}$ is invertible for every $X\in \mathcal{A}$ and every $i\in\{0,1,\cdots,k\}$.
Because $\mathsf{j}^i_{X,A} \circ q^i_{X,A}=e^i_{X,A} \circ \mathsf{i}^i_{X,A}$ we have that $e^i_{X,A} \circ \mathsf{i}^i_{X,A}$ is invertible, and since by Lemma \ref{3.1}, $\mathsf{i}^i_{X,A}$ is invertible, $e^i_{X,A}$ is also invertible for every $i\in\{0,1,\cdots,k\}$.
Thus by Proposition \ref{2.3} and Proposition \ref{3.3} $A\in (\overrightarrow{\mathcal{C}})^{\bot_{k+1}}=\mathcal{C}^{\bot_{k+1}}$.

We state the dual of Theorem \ref{1.1}.

\begin{theorem}
Let $\mathcal{A}$ be a essentially small abelian category, $\mathcal{C}$ a Serre subcategory of $\mathcal{A}$ and $q:\mathcal{A}\rightarrow \mathcal{A}/\mathcal{C}$ be the quotient functor. For an object $A\in \mathcal{A}$ and a non-negative integer $k$ the following statements are equivalent.
\begin{itemize}
\item[$(i)$]
$A\in ^{{\bot}_{k+1}}\mathcal{C}$.
\item[$(ii)$]
The natural map $q_{A,X}^i:\Ext_{\mathcal{A}}^i(A,X)\rightarrow\Ext_{\mathcal{A}/\mathcal{C}}^i(q(A),q(X))$ is invertible, for every $X\in \mathcal{A}$ and every $i\in\{0,1,\cdots,k\}$.
\end{itemize}
\begin{proof}
Apply Theorem \ref{1.1} to the opposite category $\mathcal{A}^{op}$.
\end{proof}
\end{theorem}

\section{Applications}
Let $n$ be a positive integer and $ \mathcal{M}$ be a essentially small $n$-abelian category in the sense of \cite{J}. Denote by $\eff(\mathcal{M})$ the full subcategory of $\modd\mathcal{M}$ consist of all effaceable functors, i.e. functors $F\in \modd\mathcal{M}$ with a projective presentation
\begin{equation}\label{PP}
\Hom_{\mathcal{M}}(-,Y^n)\overset{(-,f)}{\longrightarrow} \Hom_{\mathcal{M}}(-,Y^{n+1})\longrightarrow F\rightarrow 0 
\end{equation}
for some epimorphism $f:Y^n\rightarrow Y^{n+1}$. $\eff(\mathcal{M})$ is a Serre subcategory of $\modd\mathcal{M}$ (see \cite[Proposition 3.5]{EN} or \cite[Proposition 4.3]{Kv}). The composition of functors
\begin{align}
H: \mathcal{M}&\longrightarrow \modd\mathcal{M}\longrightarrow (\modd\mathcal{M})/(\eff(\mathcal{M})) \notag \\
X&\longmapsto \Hom_{\mathcal{M}}(-,X)\longmapsto \Hom_{\mathcal{M}}(-,X) \notag
\end{align}
is fully faithful \cite[Proposition 4.6]{Kv}, and it's essential image is an $n$-cluster tilting subcategory of $(\modd\mathcal{M})/(\eff(\mathcal{M}))$, see \cite[Theorem 4.7]{EN} or \cite[Theorem 7.3 and Proposition 7.5]{Kv}.
Because $(\modd\mathcal{M})/(\eff(\mathcal{M}))$ has no injective or projective object in general, proving that the essential image is $n$-rigid (i.e. for every $i\in \{1,\cdots,n-1\}$ and every $X,Y\in \mathcal{M}$ we have $\Ext^i_{(\modd\mathcal{M})/(\eff(\mathcal{M}))}(H_X,H_Y)=0$) is not easy, see \cite[Section 6]{Kv}.
In the sequel we give a different proof using Theorem \ref{1.1}.

\begin{proposition}\label{4.1}
For every $X\in \mathcal{M}$, $H_X\in \eff(\mathcal{M})^{\bot_n}$.
\begin{proof}
Let $F\in\eff(\mathcal{M})$ and consider the projective presentation \eqref{PP}. Because $f$ is an epimorphism, by the axioms of $n$-abelian categories \cite[Definition 3.1]{J}, $f$ sits into an $n$-exact sequence 
\begin{center}
$Y^0\rightarrow Y^1\rightarrow \cdots\rightarrow Y^n\overset{f}{\rightarrow} Y^{n+1}.$
\end{center}
By the definition of $n$-exact sequences we have the exact sequence
\begin{center}
$0\rightarrow H_{Y^0}\rightarrow H_{Y^1}\rightarrow \cdots \rightarrow H_{Y^n}\rightarrow H_{Y^{n+1}}\rightarrow F\rightarrow 0$
\end{center}
which is a projective resolution for $F$.
Now let $X\in \mathcal{M}$. Applying $\Hom_{\rm{mod}\text{-}\mathcal{M}}(-,H_X)$ to the projective resolution of $F$, we get the complex
\begin{align*}
0\rightarrow\Hom(H_{Y^{n+1}},H_X)\rightarrow\Hom(H_{Y^{n}},H_X)\rightarrow\cdots\\
\rightarrow \Hom(H_{Y^{1}},H_X)\rightarrow \Hom(H_{Y^{0}},H_X)\rightarrow 0
\end{align*}
which is by Yoneda lemma isomorphic to the complex
\begin{align*}
0\rightarrow\Hom(Y^{n+1},X)\rightarrow\Hom(Y^{n},X)\rightarrow\cdots
\rightarrow \Hom(Y^{1},X)\rightarrow \Hom(Y^{0},X).
\end{align*}
And by the definition of $n$-exact sequences this is an exact sequence of abelian groups. Thus we see that for every $i\in \{0,1,\cdots,n\}$ we have $\Ext^i_{(\modd\mathcal{M})/(\eff(\mathcal{M}))}(F,H_X)=0$. Because $F$ was an arbitrary effaceable functor, $H_X\in \eff(\mathcal{M})^{\bot_n}$. 
\end{proof}
\end{proposition}

\begin{corollary}
The essential image of $H: \mathcal{M}\longrightarrow  (\modd\mathcal{M})/(\eff(\mathcal{M}))$ is an $n$-rigid subcategory.
\begin{proof}
Let $X,Y\in \mathcal{M}$. By Proposition \ref{4.1} and Theorem \ref{1.1}, for every $i\in \{1,\cdots,n-1\}$ we have
\begin{center}
 $\Ext^i_{(\modd\mathcal{M})/(\eff(\mathcal{M}))}(H_X,H_Y)\cong \Ext^i_{\modd\mathcal{M}}(H_X,H_Y)=0$
\end{center}
because $H_X=\Hom_{\mathcal{M}}(-,X)$ is a projective object in $\modd\mathcal{M}$.
\end{proof}
\end{corollary}
\section*{acknowledgements}
The author wishes to express his sincere thanks to the anonymous referee for
her/his detailed comments on the previous version of this article.
This research was in part supported by a grant from IPM (No. 1400180047).


\begin{thebibliography}{10}

\bibitem{Au} \textsc{M. Auslander}, Coherent functors, In \emph{Proceedings Conference Categorical Algebra, 
} La Jolla, CA, 1965, Springer, New York, 1966, pp. 189--231.

\bibitem{CB} \textsc{W. Crawley-Boevey}, Locally finitely presented additive categories, \emph{Comm. Algebra}, \textbf{22} (1994), 1644--1674.

\bibitem{EN} \textsc{R. Ebrahimi, A. Nasr-Isfahani}, Higher Auslander's formula, Int. Math. Res. Not. IMRN. https://doi.org/10.1093/imrn/rnab219.

\bibitem{F} \textsc{C. Faith}, \emph{Algebra: Rings, Modules and Categories I}, Springer-Verlag \textbf{190}, 1973.

\bibitem{Fr} \textsc{P. Freyd}, \emph{Abelian Categories}, Harper and Row, New York, 1964.

\bibitem{G} \textsc{P. Gabriel}, Des cat\'egories ab\'eliennes, \emph{Bull. Soc. Math. France}, \textbf{90} (1962), 323--448.

\bibitem{GL} \textsc{W. Geigle and H. Lenzing}, Perpendicular categories with applications to representations and sheaves , \emph{J. Algebra}, \textbf{144} (1991), 273--343.

\bibitem{J} \textsc{G. Jasso}, n-Abelian and n-exact categories, \emph{Math. Z.}, (2016), 1--57.

\bibitem{Kl} \textsc{B. Keller}, Chain complexes and stable categories, \emph{Manuscripta Math.}, \textbf{67}(4) (1990), 379--417.

\bibitem{Kr15} \textsc{H. Krause}, Deriving Auslander's formula, \emph{Doc. Math.}, \textbf{20} (2015), 669--688.

\bibitem{KrB} \textsc{H. Krause}, \emph{Homological Theory of Representations}, Cambridge studies in advanced mathematics \textbf{195}, Cambridge University Press, 2021.

\bibitem{Kr97} \textsc{H. Krause}, The spectrum of a locally coherent category , \emph{J. Pure Appl. Algebra}, \textbf{114}(3) (1997), 259--271.

\bibitem{Kr01} \textsc{H. Krause}, The spectrum of a module category, \emph{Mem. Amer. Math. Soc.}, \textbf{149}(707) (2001).

\bibitem{Kv} \textsc{S. Kvamme}, Axiomatizing subcategories of abelian categories, \emph{J. Pure Appl. Algebra}, \textbf{226}(4) (2022), 106862.

\bibitem{OR} \textsc{U. Oberst and H. Rohrl}, Flat and coherent functor, \emph{J. Algebra}, \textbf{14} (1970), 91--105.

\bibitem{Mi68} \textsc{B. Mitchell}, \emph{On the dimension of objects and categories II. Finite ordered sets}, \emph{J. Algebra}, \textbf{9}(3) (1968), 341--368.

\bibitem{Mi} \textsc{B. Mitchell}, \emph{Theory of categories}, Academic Press, New York, 1965.

\bibitem{Po} \textsc{N. Popescu}, \emph{Abelian categories with applications to rings and modules}, Academic Press, New York, 1973.

\bibitem{Ps} \textsc{C. Psaroudakis}, Homological theory of recollements of abelian categories, \emph{J. Algebra}, \textbf{398} (2014), 63--110.

\bibitem{Ro} \textsc{J.-E. Roos}, Locally noetherian categories, In: \emph{category Theory, Homology Theory and their Applications II,} Springer Lecture Notes in Math. \textbf{92} (1969) 197--277.

\bibitem{We} \textsc{C. A. Weibel}, \emph{An introduction to homological algebra}, Cambridge studies in advanced mathematics \textbf{38}, Cambridge University Press, 1994.

\end{thebibliography}
\end{document}